# Infinitesimals, Imaginaries, Ideals, and Fictions

**David Sherry and Mikhail G. Katz**

**Abstract.** Leibniz entertained various conceptions of infinitesimals, considering them sometimes as ideal things and other times as fictions. But in both cases, he compares infinitesimals favorably to imaginary roots. We agree with the majority of commentators that Leibniz's infinitesimals are fictions rather than ideal things. However, we dispute their opinion that Leibniz's infinitesimals are best understood as logical fictions, eliminable by paraphrase. This so-called syncategorematic conception of infinitesimals is present in Leibniz's texts, but there is an alternative, formalist account of infinitesimals there too. We argue that the formalist account makes better sense of the analogy with imaginary roots and fits better with Leibniz's deepest philosophical convictions. The formalist conception supports the claim of Robinson and others that the philosophical foundations of nonstandard analysis and Leibniz's calculus are cut from the same cloth.

**§0. Introduction.** Devotees of Leibniz's calculus, e.g., the Bernoullis, were unabashedly realist about infinitesimals. Leibniz himself held a more subtle position. As Robinson puts it,

> While approving of the introduction of infinitely small and infinitely large quantities, Leibniz did not consider them as real, like the ordinary 'real' numbers, but thought of them as *ideal or fictitious*, rather like the imaginary numbers. (Robinson 1967, 33, our italics)

Robinson's observation will give pause to readers for whom "ideal" and "fictitious" are not synonymous terms. This is not to suggest that Robinson is guilty of misrepresentation. Leibniz himself writes,

> To tell the truth, I'm not myself persuaded that it's necessary to consider our infinities and infinitesimals as something other than ideal things (*choses ideales*) or well-founded fictions (*fictions bien fondées*). (Leibniz to Varignon 20 Jun 1702; GM IV, 110)

But in an earlier letter to Varignon, in which Leibniz characterizes the infinite as "ideal" and "abstract," we find,

> Yet we must not imagine that this explanation debases the science of the infinite and reduces it to fictions, for there always remains a "syncategorematic" infinite, as the Scholastics say. (2 February 1702; GM IV, 93-4; translation from Leibniz 1969, 543)

At the very least, Leibniz seems to have contemplated different ontological conceptions of infinitesimals, and he need not have settled upon a single one. We try to identify Leibniz's strongest conception. In §3 we argue that Leibniz's infinitesimals are better understood as fictions than ideal things, on the grounds that ideal things constitute possibilities, while infinitesimals do not. In §4 we dispute the syncategorematic account, voiced by a number of commentators, that Leibniz's "fictionalisms [sic] can happily be styled Archimedean" (e.g., Levey 2008, 133). Leibniz offers a conception of mathematical fictions that does not fit this description, an alternative that might be styled formalist. We argue that this conception satisfies Leibniz's deepest philosophical convictions better than the logical fictions of the syncategorematic account.

Our conclusion bears on the relation between Leibniz's calculus and Robinson's non-standard analysis. Many historians deny significant continuity between infinitesimal calculus of the 17th century and non-standard analysis of the 20th century. They point out that the latter requires resources of modern logic that were unavailable to Leibniz and his followers, and they claim

> [T]here is . . . no evidence that Leibniz anticipated the techniques (much less the theoretical underpinnings) of modern non-standard analysis (Earman 1975, 250).

and,

> The relevance of current accounts of the infinitesimal to issues in the seventeenth and eighteenth centuries is rather minimal, but it does show that it is possible to develop a consistent theory of infinitesimal magnitudes. (Jesseph 1993, 131)

Robinson sees the matter differently.

> Leibniz' approach is akin to Hilbert's original formalism, for Leibniz, like Hilbert, regarded infinitary entities as ideal, or fictitious, additions to concrete mathematics.
> (Robinson 1967, 39-40).

Like Hilbert, Robinson is a formalist, and he regards formalism as an underlying theme of his mathematics (Robinson 1970, 45; 1966, 282). §5 argues that Robinson owes a substantial debt to Leibniz, even though Leibniz did not anticipate the model-theoretic techniques of Robinson's nonstandard analysis.[1]

**§1. Imaginary Roots.** Leibniz frequently uses imaginary roots (or quantities; cf. 1899, 522) as an analogy for understanding infinitesimals.[2] For instance, in the 2 February 1702 letter to Varignon, he writes

> … even if someone refuses to admit infinite and infinitesimal lines in a rigorous metaphysical sense and as real things, he can still use them with confidence as ideal concepts (*notions ideales*) which shorten his reasoning, similar to what we call imaginary roots in the ordinary algebra, for example, $\sqrt{-2}$. (GM IV, 92; translation from Leibniz 1969, 543)

In this passage imaginary roots and infinitesimals are classified as ideal concepts (*notions*). We argue in §3.1 that ideal concepts should not be conflated with the ideal entities that Leibniz regards as constitutive of the proper subject matter of mathematics. But at this stage we are concerned to understand only Leibniz's analogy.

Imaginary quantities cannot be used to answer the questions "how many?" and "how much?" In that sense, they have no counterparts in the physical world; hence the term "imaginary" (or "impossible"). The first publication in which imaginary roots appear is Cardan's ***Ars Magna***, from 1545 (Cardan 1993; see 219-221). He presents them as solutions to an impossible problem: Divide 10 into two parts whose product is 40 (ibid., 219-20); i.e., solve the equation $x^2 - 10x + 40 = 0$. By completing the square, and operating as though the square root operation applies to negatives, Cardan obtains $5 \pm \sqrt{-15}$. The result can be 'checked' by operating arithmetically with $\sqrt{-15}$ as though it obeys the same rules as $\sqrt{15}$; thus, $\left(5+\sqrt{-15}\right)\times\left(5-\sqrt{-15}\right)=40$.

---

[1] However, Leibniz's law of continuity and transcendental law of homogeneity *do* anticipate Robinson's transfer principle and standard part function, respectively. See Katz and Sherry 2012.

[2] It is tempting to write "imaginary number" here, but this expression rarely occurs in Leibniz. Quantity, rather than number, is the object of mathematics for Leibniz. See note 20 and §5 below.

While imaginaries appear only in these pages of ***Ars Magna***, Cardan considered further cases in which they could occur. His solution for depressed cubics of the form $x^3=px+q$,

$$(1) \qquad x=\sqrt[3]{\frac{q}{2}+\sqrt{\frac{q^2}{4}-\frac{p^3}{27}}}+\sqrt[3]{\frac{q}{2}-\sqrt{\frac{q^2}{4}-\frac{p^3}{27}}},$$

yields roots of negatives for $\frac{p^3}{27}>\frac{q^2}{4}$ (103). Cardan restricts the problem so that imaginaries don't arise, presumably because he had no idea how to extract cube roots from expressions with imaginairies. Bombelli was puzzled by one such case, $x^3=15x+4$, because he knew that the roots were all real, viz., 4, $-2+\sqrt{3}$, and $-2-\sqrt{3}$, in spite of the fact that the Cardan formula gives $\sqrt[3]{2+\sqrt{-121}}+\sqrt[3]{2-\sqrt{-121}}$ as a solution. Such cases are termed "irreducible." Bombelli unraveled the puzzle once he realized that formula (1) yields a real solution if the two terms of the solution are conjugate, i.e., if they have the form $a+b\sqrt{-1}$ and $a-b\sqrt{-1}$, so that the imaginary expressions cancel when added. He was then able to show by formal manipulation that $\sqrt[3]{2+\sqrt{-121}}=2+\sqrt{-1}$ and that $\sqrt[3]{2-\sqrt{-121}}=2-\sqrt{-1}$; thus he obtained 4 as a solution in accordance with the Cardan formula. The correctness of the result was confirmed by substituting it back in the original problem.

Leibniz studied Bombelli's ***Algebra*** and was struck by the existence of imaginary expressions for real values. In a letter to Huygens of 1673 (Leibniz 1899, 550-64), he observed

> For this is the remarkable thing, that, as calculation shows, such an imaginary quantity is only observed to enter those cubic equations that have no imaginary root, all their roots being real or possible … (1899, 552; quoted in McClenon 1923, 371)

In the same letter Leibniz also observed that Bombelli failed to understand the generality of his own method; he didn't understand that a false (i.e., negative) root could result from Cardan's formula (1899, 527; quoted in Crossley 1987, 104). In the late sixteenth century mathematicians were unsure about the status of negative roots, a situation that changed after Descartes's ***Geometrie.*** In another paper Leibniz gave a general, purely algebraic proof that Cardano's formula always determines a root (1899, 560; quoted in Crossley, 106). Unlike Bombelli's check of the correctness of particular solutions, the algebraic proof does not distinguish between positive and negative values, and so demonstrates that Cardan's formula satisfies the cubic whether or not $\frac{p^3}{27} > \frac{q^2}{4}$. By not employing specific values, the argument satisfies Leibniz's criterion for a sound demonstration, viz., "reaching its conclusion by virtue of its form" thereby protecting us "against deceptive ideas" (1989, 27). The argument also manages to unite domains, the imaginary and the real, which appeared quite distinct to Cardan and Bombelli, given their inclinations to justify algebraic manipulation geometrically. In this respect Leibniz recognized a hidden power of algebra, realizing that it was more than a short hand for what could be expressed geometrically. As impressed as Leibniz was by imaginary expression for real roots, he was more impressed by the power of algebra, and more generally, symbolic thinking, to unify disparate domains. §4 argues that Leibniz thought is dominated by this idea.

The correctness of the imaginary calculations is open to an independent check, whenever the results are real values (cf. Leibniz 1920, 150). The independent check is one source for Leibniz's claim that "imaginary roots … have a real foundation" (GM IV, 93; translation from Leibniz 1969, 544). Indeed, in the following sentence he illustrates the point by referring specifically to the results obtained in the letter to Huygens. Like imaginary roots, infinities and infinitesimals, as well as dimensions beyond three and powers whose exponents are not ordinary numbers are all "founded on realities" (ibid.) in the sense of admitting an independent check. But the provenance of symbolic operations upon imaginary objects constitutes a further sense in which they are founded on realities, as symbolic operations are founded upon analogies with operations upon real quantities.[3] These analogies are responsible for the fruitfulness of introducing imaginary objects.

In the letter to Varignon, Leibniz spells out an additional connection between the real and imaginary realms.

> Even though these are called imaginary, they continue to be useful and even necessary in expressing real magnitudes analytically. For example, it is impossible to express the analytic

[3] Laugwitz understands the well-foundedness of Leibniz's fictions in this manner. "[Well-foundedness] is, of course, made precise in a principle of transfer: the *rules* are founded on the rules of "finite totalities" (1992, 152).

> value of a straight line necessary to trisect a given angle without the aid of imaginaries.[4] Just so it is impossible to establish our calculus of transcendent curves without using differences which are on the point of vanishing, and at last taking the incomparably small in place of the quantity to which we can assign smaller values to infinity. (ibid., 92)

That is, ideal concepts are necessary for the progress of mathematics. That progress takes a particular form when imaginary or ideal objects are concerned. While roots of negatives and infinitesimals appear in calculations, they are generally absent from the final result. In this sense calculations or operations with imaginary objects constitute paths, even shortcuts, from one region in the real domain to another.

**§2. Varieties of Fictions.** Costabel is uncomfortable with the analogy Leibniz draws between infinitesimals and imaginary roots (1988).[5] While the suggestion is certainly there in the 2 February letter to Varignon, Costabel argues that Leibniz remained troubled by disanalogies between the two cases and reformulated his position within the span of a few months (1988, 176-7). Leibniz insists in the 2 February letter that the science of the infinite doesn't reduce to fictions, but by 20 June of the same year he writes to Varignon that infinities and infinitesimals are fictions after all, though well-founded ones. Leibniz came to see an important distinction between imaginaries and infinitesimals.

> Le mot <u>fiction</u> qui apparaît dans le texte du 2 février 1702 est appelé par la nécessité de bien exprimer la différence entre le calcul algébrique et le calcul différentiel à travers ce qu'ils paraissent avoir en commun. Dans les deux cas, avec l'imaginaire et les infinis, on <u>nomme</u>, donc on <u>feint</u> des êtres mathématiques et il y a analogie bien exprimée par le mot <u>fiction</u>. Mais dans le second cas, celui des entitités du Calcul différentiel et integral, il n'y a pas de <u>réduction</u> à une fiction qui ne serait simplement justifiée qu'à posteriori. Il ne s'agit pas d'une <u>pure fiction</u>. (1988, 177)

---

[4] Girard discovered the analytic solution to the trisection problem, viz., $q=3x-x^3$, where q is the length of the chord that subtends the given angle (Girard 1629). The solution requires Cardan's formula. Leibniz refers to Girard's formula in the letter to Huygens (1899, 552; cf. McClenon 1923, 371).

[5] Costabel 1988 has drawn little attention. His paper, however, is the only attempt of which we are aware to come to terms with Leibniz's frequent analogies between infinitesimals and imaginaries, a project we believe to be crucial to understanding Leibnizian infinitesimals.

Leibniz, says Costabel, sees fiction as a genus, with species *pure* and *well-founded*. Both sorts of fiction involve pretending that certain mathematical things exist. And although both are "fondé en réalités," their foundations are different in kind. Imaginaries are pure fictions insofar as their justification is a posteriori. Infinitesimals are well-founded because their foundation is constituted by an identity between results obtained by their means and results "obtenus par les méthodes anciennes d'exhaustion" (1988, 176).

Costabel's analysis invites a question about fictions that are not founded upon realities. In his correspondence with Clarke, Leibniz refers, it seems, to a number of such cases: a mere will without a motive (1956, 36) , the world's having been created some millions of years sooner (38), and a material finite universe, moving forward in an infinite empty space (63). In ***Essays on Theodicy***, Leibniz characterizes wills without motives as "des Etres de Raison sans raisonnante" (G VI, 432), unreasonable beings of reason. Unreasonableness here consists in a conflict with fundamental principles; the identity of indiscernibles, for instance, "puts an end … to hundreds of fictions which have arisen from the incompleteness of philosophers' notions" (1981, 57). This list includes the tabula rasa, a vacuum, atoms, and absolute rest (ibid., 109-110). It is instructive to see how atoms conflict with this principle.

> That is why the notion of atoms is chimerical and arises only from men's incomplete conceptions. For if there were atoms, i.e., perfectly hard and perfectly unalterable bodies which were incapable of internal change and could differ from one another only in size and shape, it is obvious that since they could have the same size and shape they would then be indistinguishable in themselves and discernible only by means of external denominations with no internal foundation; which is contrary to the greatest principles of reason. (1981, 230-1)

For now, let us call fictions of the third sort *unreasonable*. Tentatively, then, there are three sorts of fiction: pure, well-founded, and unreasonable. If, as seems plausible, some fictions are *not* well-founded, then a satisfactory interpretation of Leibnizian infinitesimals should enable us to distinguish them clearly from other species of fictions.

Costabel is concerned to understand well-founded fictions, especially the sense in which they are well-founded. To aid his analysis, he considers Leibniz's comments on fictions as they first appeared in his writing, in the context of rights and politics (179-80). In that context, fictions are useful presumptions for which there is no evidence to the contrary, like presuming that a child's father is the husband of its mother (cf. Leibniz 1981, 249). This leads Costabel to identify well-founded fictions with possible ideas,

> Je crois que cette évolution de la réflexion philosophique de Leibniz aurait pu l'amener à identifier la "fiction bien fondée" avec l'"idée possible,"

and to conclude, surprisingly, that their justification is inseparable from a posteriori consequences:

> Bien fonder est resté cependant pour lui un acte qui échappe à une évidence rationnelle entière et dont la justification n'est pas séparable de conséquences à posteriori. (180)

Accordingly, well-founded fictions are creations from and in accordance with the content of experience, i.e., ideal things that have stood the test of experience. Costabel concludes, then, by denying a distinction between ideal things and well-founded fictions, and blurring his earlier distinction between pure and well-founded fictions, by putting both in the realm of the a posteriori (cf. 177). If there is something to be said for the distinction of pure from well-founded, then Costabel's analysis of infinitesimals is unsatisfactory.

**§3. Logical Fictions.** Ishiguro finds, ultimately, a clear distinction between ideal entities and fictions, and she concludes that Leibniz's infinitesimals are fictions rather than ideal entities (1991, 100). But fictionalism doesn't stand in the way of Leibniz's providing "a rigorous conceptual foundation" for infinitesimals (79-81; cf. 96). This leads her to suggest that Leibniz is a proto-Cauchy or proto-Hilbert rather than a proto-Robinson, a suggestion to which we shall return (§5).

Leibniz, according to Ishiguro, treats "infinitesimal" and "infinitely small" as non-designating expressions.

> The word 'infinitesimal' does not designate a special kind of magnitude. It does not designate at all. (83)

But this does not mean that he failed to speak significantly when he used them. Leibniz's infinitesimals are fictions, according to Ishiguro, because "infinitesimal" fails to refer; but they are *well-founded* fictions, and so meaningful even in scientific discourse. Ishiguro takes pains to specify the sense of "well-founded."

> … although the limit and the differential may be an ideal entity or a well-founded fiction, the language of infinitesimals or the language of limit is not, according to Leibniz, a fiction. … [E]very proposition or equation which makes reference to these ideal entities or fictions can be paraphrased into propositions and equations in which such entities are not designated, which have strict truth conditions. (92-3)

There are, in other words, no infinitesimal magnitudes, but propositions apparently referring to such things can state truths since they can be paraphrased as propositions that refer only to finite magnitudes. Thus there are strict truth conditions accompanying the paraphrases.

Ishiguro's interpretive strategy employs what Russell called logical or symbolic fictions (1919, 45 and 184). In fact, Ishiguro sees in Leibniz an anticipation of Russell's theory of descriptions:

> Leibniz took a position closer to that of Bertrand Russell in his theory of descriptions, a theory which Russell claimed was a theory of contextual definition. … [Leibniz] gives a contextual definition to the word "infinitesimal" as well as to the sign for differential quotient. (82)

Russell's theory of descriptions explains how a proposition containing an apparently non-referring description, "the present king of France is bald," has definite truth conditions in spite of the non-referring expression. The truth conditions emerge once the proposition's *true* logical form is made explicit. The apparent logical form is a property being affirmed of an object, viz., x is F; but the true logical form is an existentially quantified formula. Thus the proposition is false because there is no object – among the ordinary objects of perception - satisfying the formula.[6] Russell applies the same gambit to propositions that refer, only apparently, to classes (1919, ch. 17). Once the true logical form has been made explicit, we see that such propositions refer to ordinary objects. Leibniz is making precisely the same move, according to Ishiguro.

---

[6] To wit: "x is king of France & $\forall y$(y is king of France$\supset$y=x) & x is bald." We ignore the complications faced by Russell's theory in view of "the present king of France is not bald."

She supports this reading with various passages, including one from the 2 February letter to Varignon:

> … we must consider that these incomparable magnitudes themselves, as commonly understood, are not at all fixed or determined but can be taken to be as small as we wish in our geometrical reasoning and so have the effect of the infinitely small in the rigorous sense. If any opponent tries to contradict this proposition, it follows from our calculus that the error will be less than any possible assignable error, since it is in our power to make this incomparably small magnitude small enough for this purpose, inasmuch as we can always take a magnitude as small as we wish. (1969, 543; cf. 1966, 28)

An infinitesimal expression, "dx" for example, does not refer to anything. However, a proposition in which it occurs is still significant.

> Leibniz is saying that whatever small magnitude an opponent may present, one can assert the existence of a smaller magnitude. In other words, we can paraphrase the proposition with a universal proposition with an embedded existential claim. This fact about the existence of smaller and smaller finite magnitudes is supposed to justify certain rules governing the sign 'dx'. (Ishiguro 1990, 87)

Thus, the equation that gives the slope of the tangent to $y = \frac{x^2}{a}$, viz., $\frac{dy}{dx} = \frac{2x}{a}$, refers, apparently, to a ratio between infinitesimal magnitudes. This equation can be paraphrased as a universal with an embedded existential claim:

> <u>for any</u> finite error ε you may assign, <u>there is</u> a finite difference dx in the x variable and <u>there is</u> a (corresponding) finite difference dy in y variable, such that the difference between $\frac{dy}{dx}$ and $\frac{2x}{a}$ is less than ε.

Thanks to the paraphrase, we see that the equation designates only finite magnitudes, ε as well as dx and dy, whose choice depends upon ε. The ratio $\frac{dy}{dx}$, then, is a logical fiction. It is a fiction because it does not designate an entity, and it is a logical fiction because the real content of the propositions in which it occurs is disclosed by the proposition's true logical form.[7] This analysis of infinitesimals – the syncategorematic interpretation of infinitesimals - is a starting point for much recent Leibniz scholarship (e.g., Arthur 2008, 20 and Levey 2008, 107).

Although Ishiguro's discussion of infinitesimals is originally non-committal on the issue of whether they are ideal entities or fictions (cf. 92-3), by the end of the chapter she draws a firm distinction between those categories and places infinitesimals in the latter.

> One must distinguish what Leibniz called ideal entities, *entia rationis* or *êtres de raison*, on the one hand and fictions, *entia ficitium*, on the other. When we refer to numbers or to relations which are ideal entities, we are referring to something. … Fictions, on the other hand, are not entities to which we refer. They are not abstract entities. They are correlates of ways of speaking which can be reduced to talk about more standard kinds of entities. Reference to mathematical fictions can be paraphrased into talk about standard mathematical entities. (100)

Ishiguro's paradigm ideal entities are relations, finite numbers, and magnitudes (99-100). Unlike their infinite and infinitesimal counterparts, finite numbers and magnitudes are "specifiable" (French *assignable*) (cf. G VI, 90, quoted in Ishiguro, 79); as Ishiguro explains, they "have a secondary existence dependent upon the possible existence of substances" (99) – they are "properties of possible individuals" (100). Expressions for fictional entities have no referents because such referents are not possible; if they were possible they would, like (ordinary) numbers, be properties of possible objects.

**§3.1. Ideal Entities.** This interpretation puts Ishiguro at odds with Costabel, for whom well-founded fictions are possibilities. Much evidence favors her interpretation. In a letter to Bernoulli Leibniz writes

> … as concerns infinitesimal terms, it seems not only that we cannot penetrate to them, but

---

[7] Ishiguro offers a tidy account of Leibniz's conception of infinitesimals. Its tidiness results, however, from a logical theory that was neither available to Leibniz nor amenable to his own logical theory. Russell's logic is a logic of relations, while Leibniz, following Aristotle, uses a subject/predicate logic. Moreover, Russell's logic is extensional while Leibniz's was, by many accounts, intensional. We ignore these complications, presuming that Leibniz could paraphrase the Russellian idiom into his own.

> that there are none in nature, that is, that they are not possible. (G.M. III, 499-500; quoted in Ishiguro, 84)

The impossibility Leibniz envisions here is not explicit, but elsewhere he suggests that if he were to accept infinitesimals, he would also be committed to infinite lines terminated on both ends (1969, 543). They would be inverses of infinitesimals, and the impossibility of an infinite but terminated line is, perhaps, more obvious. Ishiguro also treats the limit "regarded as being included in the value of the variables that it limits" as a well-founded fiction[8] (90), and here too the fiction looks to be impossible for Leibniz. For he writes, "it is not at all rigorously true that … a circle is a kind of regular polygon" (Leibniz 1969, 546).

We are in substantial agreement with Ishiguro's depiction of ideal entities as possibilities and mathematical fictions as impossibilities. But it requires elaboration in light of some apparently conflicting remarks of Leibniz. Compare Leibniz's view that infinitesimals are impossible because there are none in nature with a similar remark about continuity in the 2 February letter to Varignon

> … one can say in general that, though continuity is something ideal and there is never anything in nature with perfectly uniform parts, the real, in turn, never ceases to be governed perfectly by the ideal and the abstract … . (GM IV, 93; translated in 1969, 544)

This is troubling because continuity is *not* one of Leibniz's stock examples of a fiction. But if failure to appear in nature is a sign of impossibility, then a continuous magnitude, which requires uniform parts[9], ought to be as fictional as an infinitesimal. Moreover, the virtue of continuity, its capacity to govern the real perfectly, is shared by infinitesimals, as Leibniz notes earlier in the same letter.

> So it can also be said that infinites and infinitesimals are grounded in such a way that everything in geometry, and even in nature, takes place as if they were perfect realities. (ibid.)

Leibniz treats imaginary roots the same way in an earlier letter to Bernoulli.

> For perhaps the infinite, such as we conceive it, and the infinitely small, are imaginary, and yet apt for determining real things, just as imaginary roots are customarily supposed to be.

[8] The limit itself, e.g., a circle that limits a series of regular polygons with increasingly more sides, is not a fiction; rather it is an ideal entity.
[9] Elsewhere Leibniz describes the continuous as a repetition and diffusion of the same nature (1989, 251).

> These things are among the ideal reasons by which, as it were, things are ruled, although they are not in the parts of matter. (June 1698 letter to Bernoulli, GM III 499-500; translated in Jesseph 1998, 28)

These are not the only passages in which Leibniz appears to blur the distinction between ideal and fictional. After dismissing unreasonable fictions (e.g., the world's having been created millions of years sooner) as "altogether unreasonable and impracticable," and insisting that they "cannot be admitted," Leibniz asserts

> Mere mathematicians, who are only taken up with the conceits of imagination, are apt to forge such notions; but they are destroyed by superior reasons. (1956, 63-4)

Yet elsewhere Leibniz is less dismissive of the notions that mathematicians have forged. Commenting in ***New Essays*** upon "hundreds of other fictions which have arisen from the incompleteness of philosophers' notions" he is more conciliatory.

> They are something which the nature of things does not allow of. They escape challenge because of our ignorance and our neglect of the insensible; but nothing could make them acceptable, short of their being confined to abstractions of the mind … Whereas abstraction is not an error as long as one knows what one is pretending not to notice is *there*. This is what mathematicians are doing when they ask us to consider perfect lines and uniform motions and other regular effects, although matter (i.e. the jumble of effects of the surrounding infinity) always provides some exception. This is done so as to separate one circumstance from another and, as far as we can, to trace effects back to their causes and to foresee some of their results. (1981, 57)

Here the forged notions of mathematicians are the abstract or ideal objects of mathematics instead of inadmissible conceits. Such passages threaten the distinction between ideal entities and fictions and so prompt a closer look at possibility and impossibility for mathematical concepts, and their failure to occur in nature.

It is tempting to understand possibility in light of Leibniz's idea of analysis. A possible notion is one whose analysis leads to notions known to be consistent and compatible with one another (e.g., 1989, 26). In the case of adequate knowledge, the analysis is carried to completion, i.e., reduced to "*primitive possibilities* or to irresolvable notions" (ibid.). Leibniz says he "won't now venture to determine" whether people can ever carry an analysis to completion. But he does allow *additional* means for demonstrating possibility, viz., a causal definition whereby "we understand the way in which a thing can be produced" (ibid.; cf. 1969, 230). This is a traditional notion of mathematical possibility (or existence[10]), although the concept is not completely determinate, owing to disagreement over the suitable means of production. The causal definition seems to be the practical benchmark of possibility for Leibniz. While adopting a wider notion of suitable means of production than, say, Euclid or Descartes[11], Leibniz still requires some sort of production in the imagination. This criterion rules out infinitistic concepts like the largest number or the smallest fraction, which would require an infinite number of steps to produce.

The production of mathematical objects requires a material or content; for Leibniz, the material is sensible. Indeed, all thought, even abstract thought, requires sensation, according to Leibniz (G IV, 563; cf. McRae 1995, 179). Some sensible qualities, like number and shape, are shared by more than one sense, and our concepts of shared sensible qualities are found in the internal sense or imagination (1969, 548). Unlike our concepts of colors, tastes, etc., concepts of number, shape, etc., can be known distinctly, because they admit definition.[12] Hence, these concepts constitute the subject matter of the mathematical sciences (ibid.). They are, as noted previously, ideal and abstract concepts, corresponding to nothing in nature, because nature does not admit the uniformity which characterizes mathematical concepts; the appearance of uniformity, i.e., continuity, arises only from deficiency in our senses (1981, 110; cf. McRae 1995, 182). Thus the subject matter of mathematics consists in imaginable but ideal objects, which arise from sensation but require a contribution from the understanding as well (1969, 548). A logically coherent mathematical concept, then, is one definable as a construction from sensible content. That is, a proper mathematical concept has sensible content, and so, imaginary concepts are not properly concepts at all (even though imaginable concepts *are*).[13] This conclusion makes sense of Leibniz's disparaging descriptions of ordinary mathematical concepts as "monstrous" and "something which the nature of things does not allow of" (G II, 249; 1981, 57). Ideal and abstract concepts are only approximations of nature's infinite variety, and in that sense they are fictional. They can occur in nature as it appears to us, but they can't occur in nature as it is; if we were aware of our petite perceptions, and could accommodate them in our theories, we would avoid incomplete notions. But genuine mathematical concepts – ideal entities - must arise from sensation, and without their sensible provenance they could not be objects of mathematical science, for there would be nothing to which they refer. In that sense ideal entities are different from mathematical fictions.

---

[10] Leibniz states that if infinitesimals were possible, then he would be forced to grant their existence (GM III, 524; quoted in Ross 1990, 128).

[11] Euclidean construction is confined to ruler and compass, while Cartesian construction includes any algebraic magnitude. Leibniz extended that domain to any magnitude describable by a precise and continual movement. See Breger 2008, 149.

[12] In ***New Essays*** Leibniz observes that the geometries of the blind man and the paralytic "must come together, and agree, and indeed ultimately rest on the same ideas, even though they have no images in common" (137). Thus, mathematical concepts require more than sensation for their possibility.

Infinitesimals and imaginary roots are not, therefore, *objects* of mathematical science, and their fictional status rests ultimately upon this absence of objectivity. These concepts, because they contain a contradiction, prevent us from imagining objects in accordance with their definitions. But even though thought requires sensation, our inability to imagine infinitesimal quantities or imaginary roots does not prevent our reasoning symbolically with them (1989, 25). Symbols provide enough sensible content to make reasoning possible; that is, the symbol enables us to think of an object without having an idea of it (ibid., 25-6). And even when we have an idea, for the most part our thinking about complex matters, like mathematics, is symbolic. But at the core of symbolic thought there lies a sensible content, some of whose properties can be revealed by manipulation of a purely symbolic structure.

Ishiguro's distinction between ideal and fictional entities remains intact, then. While continuous magnitude, for instance, is impossible with respect to nature, it is nonetheless possible with respect to our sensible representation of nature; for we are able to construct continuous magnitude from the content of sensation. Imaginary roots and infinitesimals are a different story, however; because they contain contradictions, it is impossible to construct corresponding objects in imagination. Mathematical propositions, therefore, contain two sorts of expressions: those that refer to the content of sensation and those that don't, at least not directly. Clearly the ideal concepts (*notions ideales*) of the 2 February 1702 letter to Varignon should be treated as expressions of the *second* type. Whether or such fail to refer even *in*directly is the concern of §4.

---

[13] Ross argues that the reality of Leibniz's 'real' mathematical concepts is not a function of their being instantiated in physical reality. Geometrical figures cannot exist in reality because "physical objects are distinguished by greater or lesser deviations from perfection" (131). Odd as it may sound, incomplete notions involve a sort of perfection, as when two lines are exactly equal to one another, or two real objects have exactly the same size and shape (G II, 249). In spite of Leibniz's charge of monstrosity, incomplete notions *are* logically coherent according to Ross: "They are *entia rationis*, or "mental entities"" (133), that is, ideal entities in Ishiguro's sense.

**§3.2. Unintelligible Fictions.** Ishiguro asserts two theses about Leibnizian fictions: (1) they are all impossibilities, and (2) they are all *logical* fictions. Indeed, in explicating *entia fictium* she asserts

> Fictions … are correlates of ways of speaking which can be reduced to talk about more standard kinds of entities.

Both theses are too strong, and they conflict with her subsequent remark that the world's having been created millions of years sooner is an "impossible fiction," distinct from "useless" fictions like Pegasus (144).[14] For this distinction to be meaningful, Pegasus must be a possible fiction, and, generally speaking, literary, legal, etc. fictions seem to refer to possibilities, thanks to their similarity to actual cases. More importantly, Ishiguro says that the world's having been created millions of years sooner constitutes an "unintelligible hypothesis" (ibid.). But unintelligible propositions are *neither* true *nor* false, and so one could not paraphrase such a proposition to obtain a new proposition having 'strict truth conditions'. Owing to these difficulties it is better to limit Ishiguro's remarks on fictions to mathematical (logical, scientific) fictions: They are impossibilities, but propositions in which they occur may be paraphrased as propositions that refer to standard, i.e., possible mathematical entities.

We argued in §2 that an explication of well-founded fictions should distinguish them from other kinds, viz. pure fictions and unreasonable fictions. Ishiguro mentions instances of each, imaginary roots (99) and, as we saw, the world's having been created millions of years sooner (143-4); but she makes no attempt to paint either of them as *logical* fictions. In each case, her comments raise issues for an analysis of Leibnizian infinitesimals. First, unreasonable fictions constitute a different sort of impossibility from that of well-founded fictions; the latter are, presumably, intelligible, whereas the former are not. This poses a formidable puzzle, which we attempt to solve in §4: In what does the intelligibility of an impossibility consist? We can, however, specify what Leibniz understands by intelligibility.

[14] The suggestion that all fictions are logical may be the result of adding a chapter on infinitesimals to the second edition of the book, eleven years after the first edition.

In a passage from the Leibniz-Clarke correspondence Leibniz criticizes a number of fictions that Jesseph describes as "Newtonian monstrosities," like the void and action at a distance (1998, 35). The passage occurs in the course of defending the principle of sufficient reason.

> Has not everybody made use of this principle, upon a thousand occasions? 'Tis true, it has been neglected, out of carelessness, on many occasions: but that neglect, has been the true cause of chimeras; such as are (for instance,) an absolute real time or space, a vacuum, atoms, attraction in the scholastic sense, a physical influence of the soul over the body, [and of the body over the soul] and a thousand other fictions, either derived from erroneous opinions of the ancients, or lately invented by modern philosophers. (1956, 95-6)

This passage treats the principle of sufficient reason as a necessary condition of theoretical intelligibility, and, presuming that Leibniz does not mean to disparage *all* fictions, it provides a necessary condition of intelligibility for fictions too. A well-founded fiction, then, would have to satisfy the principle of sufficient reason as well as withstanding the test of experience. This conclusion will play heavily in our interpretation of Leibniz's fictionalism. Of course, we are still in need of a sufficient condition for an intelligible fiction.

**§3.3. Non-logical Fictions.** Ishiguro's comments on imaginary roots pose a different problem. She notes Leibniz's analogy between infinitesimals and imaginary roots (84-5; 99), and her analysis of *entia fictium* suggests that, like the former, the latter *are* eliminable by paraphrase. Yet she neither finds one in Leibniz nor offers one of her own.[15] This is unsurprising, perhaps, for the reason that led Costabel to deny imaginary roots the status of well-founded fictions in the first place:

---

[15] Referring to infinitely large and infinitely small magnitudes, Leibniz writes, "For I hold both to be fictions of the mind due to an abbreviated manner of speaking, fitting for calculation, as are also imaginary roots in algebra" (G II, 305). This can be read as a claim that imaginary roots, like infinitesimals, are eliminable by paraphrase. But it may also be read as asserting that imaginary roots, like infinitesimals, are fitting for calculation. To claim that imaginary roots are no more an abbreviated manner of speaking would require an analogy to Cardan's formula, which gave solutions to the irreducible case of the cubic without the appearance of imaginaries.

Their results are not generally obtainable by methods that employ only ordinary quantities. The 'ancient' method of exhaustion uses ordinary quantities to determine quadratures, so it is not surprising that quadratures obtained by infinitesimals might be paraphrased finitistically. But imaginary roots had no such predecessor. They were invented to provide solutions where none had existed before; hence Cardan's notorious comment that imaginaries are as subtle as they are useless (1993, 220). *Prima facie*, then, there is no reason to expect that propositions containing imaginary roots should be reducible to propositions containing only ordinary quantities. This situation muddies the claim that Leibnizian infinitesimals are well-founded fictions. Given that some commentators treat imaginaries as well-founded fictions (e.g., Ross 1990, 128; Jesseph 1998, 35) and given that a thorough analysis of Leibnizian infinitesimals should specify their relation to other fictions, we are led to inquire why Leibniz prized imaginaries so highly. Once this is clear, we will be in position to appreciate why the syncategorematic interpretation of infinitesimals does not satisfy Leibniz's deepest philosophical convictions.

**§4. Pure Fictions Revisited.** If the syncategorematic interpretation of Leibniz's fictionalism is adequate, any proposition that refers to a fiction can be paraphrased as a proposition that refers only to the sensible content of mathematics. This interpretation accounts neither for Leibniz's conception of mathematical fiction, nor even for his conception of infinitesimal. Charity bids us to look for a better interpretation.

A *technical* difficulty with the syncategorematic interpretation of infinitesimals is that it doesn't deliver what it promises: Infinitesimals aren't fully eliminable. According to Jesseph, Leibniz never "attempted anything like a general proof of the eliminability of the infinitesimal, or offered anything approaching a universal scheme for re-writing the procedures of the calculus in terms of exhaustion proofs" (2008, 233; cf. also Mancosu 1996, 170). In a later discussion of ***De Quadratura Arithmetica***, the zenith of Leibniz's syncategorematic efforts, Jesseph explains why such a universal scheme was not forthcoming: ***De Quadratura*** depends on infinitesimal resources in order to construct an approximation to a given curvilinear area less than any previously specified error.

> Although this sort of procedure can show the eliminability of infinitesimals from a large class of problems, it falls short of a completely general method that could apply to any curve

> representable by an analytic equation. The problem is that the construction of the auxiliary curves requires that we have a tangent construction that will apply to the original curve. This is readily available in the case of the circle, and tangents to conic sections and other well behaved curves constructible with classical methods. However, one great strength of the infinitesimal calculus is that it offered algorithmic solutions to problems of tangency for a wide variety of curves, yet the procedure in the Arithmetical Quadrature could only be made fully general if we already had a solution to the general problem of tangent construction. (2012)[16]

This problem is reminiscent of the difficulty that led to infinitesimal methods in the first place. Archimedes' method of exhaustion required one to determine a value for the quadrature in advance of showing, by reductio argument, that any departure from that value entails a contradiction. Archimedes possessed an heuristic, indivisible method for finding such values, and the results were justified by exhaustion, but only after the fact. That is, Archimedes had no general proof of the validity of his indivisible method. By the same token, the use of infinitesimals is 'just' a shortcut only if it is entirely eliminable from quadratures, tangent constructions, etc. Jesseph shows that this is not the case. So, if the substance of Leibniz's fictionalism is logical fictionalism, in Ishiguro's sense, it is a failed doctrine.

The failure to execute the syncategorematic strategy is a technical failure, the sort that can be remedied by technical advance. In this respect it resembles a number of Leibniz's proposals, not the least of which is his plan to reduce mathematics to logic. The latter lacked a sufficiently developed logic of relations, which appeared only two centuries later, in the work of Frege and Peirce. Around the same time, Dedekind, Weierstrass, and Cantor introduced the technical advances - especially the notion of infinite aggregate – necessary to carry off the syncategorematic strategy. In §5 we consider the extent to which the reduction of analysis to arithmetic by this 'triumvirate' is consistent with Leibniz's syncategorematic strategy. For now it suffices to note that even though Leibniz did not carry off the syncategorematic strategy, it illustrates his talent for anticipating future developments.

---

[16] For a specific instance of the auxiliary construction, see Bos 1980, 63-4.

Yet technical failure is not the gravest shortcoming of the syncategorematic strategy. Simply put, it fails to account for the full range of mathematical fictions that Leibniz recognizes; in particular, it has nothing to say about the fictional quantities to which algebraic investigations lead – not only imaginary roots, but negative quantities (a mathematical fiction that doesn't draw as much attention as Leibniz's standard examples) as well (cf. "Observatio quod Rationes," GM V, 387-9).[17] There is no indication that Leibniz is prepared to treat negative or imaginary quantities as logical fictions.[18] If, following Costabel, we refer to these quantities as pure fictions, our objection is that the syncategorematic strategy fails to account for pure fictions.[19]

That negative quantity involves an absurdity is already apparent from Leibniz's Latin, *quantitas nihilo minor*, quantity less than nothing. But Arnauld, who helped to introduce Leibniz to mathematics, presented a more interesting paradox for negative quantity in questioning whether minus times minus should yield plus.

> I do not know how to fit this with the fundamental property of multiplication, that is that the unit is to one of the magnitudes that are multiplied as the other is to the product … For shall one say that +1 is to -4 as -5 is to +20? I do not see this. For +1 is greater than -4. On the contrary -5 is less than +20. (quoted in Mancosu 1996, 88)

The idea that a greater is to a lesser as a lesser to a greater was enough to convince Leibniz that these ratios are not real but imaginary (GM V, 388). He adds that he would not deny that -1 is a quantity less than nothing, provided this is understood *sano sensu* – in a proper sense (ibid.). Leibniz goes on to compare his denial of the reality of quantity less than nothing with his denial of the reality of the infinitely large and the infinitely small (389), suggesting again that there is no harm in speaking that way if it is understood *sano sensu*. The best account of Leibniz's fictionalism, in our view, is one that elaborates how fictions are *properly understood*.

---

[17] Robinson comments that Leibniz did not consider infinitesimals "as real, like the ordinary 'real' numbers" (1967, 33). It would be interesting to discover whether Robinson interprets Leibniz as excluding negatives from the domain of the mathematically real (i.e., *entia rationis*).

[18] Grosholz and Yakira present a text from Leibniz as containing "a definition of negative numbers" (1998, 84; cf. 89-98 where the Latin text, LH XXXV, 1, 9, pp. 9 recto-14 recto, is reproduced). It is an heuristic explanation of the concept, rather than a contextual or operational definition in Ishiguro's sense. The authors contend that the definition exhibits the conditions of intelligibility of negative numbers by an analogy with not only geometry but also kinematics (ibid.).

[19] While Levey is sensitive to different types of mathematical fictions in Leibniz, his conclusion that Leibniz's fictionalisms can be styled Archimedean (2008, 133) depends on ignoring both negative and imaginary quantities.

It is tempting to suppose that the proper sense of a proposition involving fictions lies in a paraphrase that eliminates the fiction, as Leibniz does nothing discourage this interpretation. However, a reduction of negative quantities to admissible entities, able to handle ratios, was unavailable until Hankel's work in the 1860's, and only after mathematicians had given up the idea that quantity was the object of mathematical science (Schubring 2005, 601-2; cf. Hankel 1867).[20] The same is true of imaginary roots, in spite of claims (e.g., Burgess and Rosen 1997, 222-3) that their geometric interpretation erased any doubts about their reality.[21] The elimination strategy was no doubt attractive to Leibniz, as he was forced to defend the infinitesimal calculus from critics who admired traditional mathematics (cf. Mancosu, 1996, 165ff.). And granted, negative and imaginary quantities never came under the same kind of scrutiny. But an account of mathematical fiction that ignores imaginary quantities is unsuitable as an interpretation of Leibniz. It ignores the fact that imaginary quantities are a paradigm mathematical fiction for Leibniz. And in doing so it lacks the universal scope that his philosophical doctrines urged. His was not a port-manteau philosophy; Leibniz sought universal harmony, and, as Knobloch puts it, "Every harmony implies generality, while generality implies beauty, conciseness, simplicity, usefulness, fecundity" (2008, 181). It is outside of Leibniz's philosophical character to suppose that mathematical fictions come in two distinct varieties. In other words, the syncategorematic interpretation gives short shrift to connections between Leibniz's philosophy and his mathematics, ignoring virtues shared by infinitesimals, imaginary roots, negative quantities, and other mathematical fictions. The proper sense of propositions that refer to fictions - Leibniz's primary account of mathematical fictions, in our view - emerges in light of these virtues.

---

[20] It is for this reason that we stress in note 2 Leibniz's use of "imaginary quantity" instead of "imaginary number." The impossibility of imaginary roots is more palpable in the case of quantity than in the case of number.

[21] Otherwise we could not make sense of the importance of Bolzano's purely analytic proof of the intermediate value theorem, which sought to rid the fundamental theorem of algebra of any dependence upon geometry (cf. Kitcher 1975). Indeed, the project of formulating arithmetic and algebra independently of geometry goes back to Monge and Lagrange in the eighteenth century (Martínez 2006, 81).

What are these virtues? The 2 February letter to Varignon observes that fictions – not just infinitesimals, but imaginary roots, powers beyond three and even powers whose exponents are not ordinary numbers - are introduced "all in order to establish ideas fitting to shorten our reasoning." (1969, 543). In a letter to Dangicourt Leibniz praises fictions for abbreviating our thought and for enabling us to speak universally (Leibniz 1768, III, 500-1, translated in Jesseph 2008, 230). And in his treatment of negative quantities and their ratios, he says that such fictions "have great utility in calculation and discovery (the art of invention) and they serve as universal concepts" (GM V, 388). Each of these slogans emphasizes convenience, and it is tempting to imagine that convenience is the primary rationale for mathematical fictions, as if the following, memorable description from the Scottish philosopher William Hamilton[22] had some merit.

> The mathematical process in the symbolical method [i.e., the algebraical] is like running a railroad through a tunneled mountain; that in the ostensive [i.e., the geometrical] like crossing the mountain on foot. The former causes us, by a short and easy transit, to our destined point, but in miasma, darkness, and torpidity, whereas the latter allows us to reach it only after time and trouble, but feasting us at each turn with glances of the earth and of the heavens, while we inhale the pleasant breeze, and gather new strength at every effort we put forth. (quoted in Olson 1971, 42-3).

Leibniz would have disagreed heartily, and while he may not have employed pastoral metaphors, he would certainly have extolled the aesthetic virtues of the algebraic method.

We introduced negative quantities and their ratios, because they are relevant to some useful insights into the virtues of mathematical fictions. Mancosu observes that in spite of its conflict[23] with the intuitive understanding of the notion of ratio, Leibniz is comfortable *extending* the fractional calculus to negative numbers as long as they are understood as imaginary quantities (1996, 91). In so doing, Leibniz "takes advantage of the algebraic mode of thought to unify what ought to be kept distinguished only at the foundational level." Leibniz's rationale is intriguing.

---

[22] This is not William *Rowan* Hamilton (1805-1865), the Irish physicist and mathematician, who discovered quaternions, but the Scottish metaphysician who lived from 1788-1856.

[23] Mancosu's text says "conflate," but we assume this is a typographical error.

> I admit that I have always calculated the analogies by fractions. For to what use should one introduce a new sort of characters, as if the analogy y.a:a.x was something new whose calculus were not subsumed in the common precepts of arithmetic or algebra. (quoted by Mancosu 1996, 91).

As Leibniz sees it, a calculus that failed to take advantage of existing calculi has little value. Mathematical progress comes, in other words, by grafting new concepts on to existing practices. A second insight helps to understand that analogy. In discussing (L'Hôpital's presentation of) Leibniz's calculus, Mancosu remarks

> With respect to the Cartesian tradition, the two postulates of L'Hôpital represented a concept-stretching[24]: he stretched both the notion of equality, considered now as a relation between two quantities that differ by an infinitely small quantity, and the notion of polygon, extended now to encompass curves. (152-3)

Concept-stretching here consists in extending - to a new context - techniques that work in a familiar one. Treating a new case as an instance of a more familiar one allows us to take advantage of intuitions that accompany the familiar case. Leibniz's use of the characteristic triangle illustrates the point nicely. Leibniz and his followers use the characteristic triangle to determine the tangent to a curve at a given point, P. It is an infinitesimal right triangle whose hypotenuse is the tangent and whose sides are the infinitesimal increases in the abscissa and ordinate, dx and dy. The characteristic triangle is similar to the triangle consisting of the tangent, the ordinate and the subtangent; thus, computing $\frac{dy}{dx}$ allows one to determine the subtangent and so construct the tangent. While the characteristic triangle is a limit entity (the limit of secants through P and points progressively closer to P) unavailable to intuition, *pretending* the tangent at P is the hypotenuse of such a triangle allows us to draw upon geometrical intuition, and so take advantage of existing mathematics in the development of further mathematics.

Concept-stretching is similarly evident in the introduction of negative quantities and imaginary roots, though, of course, these fictions were not introduced by Leibniz. Although there are no clear intuitions of negative quantities,[25] they are treated as though they obey rules for ordinary quantities, e.g., equals added to equals gives equals, etc. Of course, there are cases, like the multiplication of two negatives, to which the ordinary rules do not apply unambiguously. Rules governing those cases must be crafted in order to accomplish specific ends.[26] For Leibniz the end is universal harmony, or the greatest possible systematization, and certainly the introduction of negatives achieves that end. By embracing negative quantities, one can reduce 13 different forms of the cubic equation[27], each requiring a different strategy for solution, to a single form solvable by a single strategy. Likewise, by embracing imaginary roots, one extends the scope of Cardan's formula to cubics for which the discriminant is non-negative; that is, one no longer has to recognize exceptions. But the systematization made possible by imaginary roots is even greater, as Girard grasped when he first articulated the fundamental theorem of algebra (Girard 1629[28], translated in Struik 1969, 85-6). Imaginary roots allowed Girard to turn his strategy for finding roots of polynomials from factions (i.e., sums of products of the coefficients) into the plausible conjecture that *every* polynomial has the same number of roots as the degree of the polynomial.[29] It is difficult to imagine a simpler and more aesthetically pleasing conclusion to an episode that commenced with numerous exceptions and special cases.

[24] Concept-stretching, of course, is crucial to Lakatos's account of mathematical progress (cf. Lakatos 1976).

There is a further respect in which imaginary roots contribute to systematization: They avoid an abrupt discontinuity in favor of a transition to a different state. In the absence of imaginary roots, Cardan's formula ceases to determine roots once $\frac{p^3}{27}$ exceeds $\frac{q^2}{4}$. For a rule or concept to become abruptly inapplicable is anathema to a philosopher whose commitment to continuity of transition is so strong that he denies the possibility of generation or death in any strict sense (cf.

---

[25] There are intuitions corresponding to negative quantities, e.g., debts; but they are not clear or consistent intuitions. Thus, it is counterintuitive that multiplying two debts should yield an asset, and that the square root of a debt should be imaginary, and so neither less than nor greater than zero.

[26] Since Hankel's ***Theory of Complex Numbers*** (1867) preserving distribution has been the usual rationale for the rule that negative times negative is positive (cf. Fischbein 1987, 99). If, e.g., $-1\times-1$ were $-1$, then applying distribution to $-1\times(1-1)$ would yield -2 instead of 0.

[27] Each of these forms is considered separately in chapters 11-23 of Cardan's A***rs Magna***.

[28] Girard's text does not have usual page numbers. The passage occurs in the explication of Theorem II, under "Des Equations ordonnées."

[29] This is not to suggest that Girard could prove the fundamental theorem of algebra. But his strategy for finding roots by means of the factions was a giant step in that direction.

***Monadology*** §73, 1989, 222). Imaginary roots insure that the arithmetic operations continue to have meaning and application even when $\frac{p^3}{27}$ exceeds $\frac{q^2}{4}$. Thus imaginary roots satisfy a scientific ideal that is at least as important to Leibniz as mollifying his critics.

> I recognize that I myself attach great importance to those who endeavour to bring carefully all the demonstrations back to their first principles and to have constantly devoted to this all my efforts. But for all that I do not say to hinder the art of invention on account of too many scruples, nor to reject under this pretext the best of discoveries, by depriving ourselves of their advantages. (GM V, 322, quoted in Mancosu 1996, 90)

Even though propositions with imaginary roots are not reducible to propositions that concern just ideal objects, the contribution of imaginary roots to the systematization of algebraic thought – along with their utility in finding real roots – establishes them as too important to be abandoned for the sake of a traditional ideal of mathematical rigor. At issue here is nothing less than Leibniz's principle of continuity, which enjoins us to find underlying forms that transcend apparent differences and so promote the end of systematization.

Leibniz ordinarily invokes the principle of continuity in the context of quantities approaching one another (cf. 1969, 351-4), but he is also prepared to use the principle in wider contexts.[30]

> … what we call birth or death is only a greater and more sudden change than usual, like the drop of a river or of a waterfall. But these leaps are not absolute and not the sort I disapprove of; e.g., I do not admit a body going from one place to another without passing through a medium. Such leaps are not only precluded in motions, *but also in the whole order of things or of truths*. (Leibniz 1951, 188; authors' emphasis)

Both negative quantities and imaginary roots move us in the direction demanded by the principle of continuity, a principle that mirrors the best of all possible worlds: "the way of obtaining as much variety as possible, but with the greatest order possible" (***Monadology*** §58, 1989, 220).

---

[30] "All natural connections, whether causal, spatiotemporal, or conceptual, are governed by this principle" (Jorgenson 2009, 223). In illustrating the law of continuity in "Cum prodiisset," Leibniz appeals to the analogy between points at infinity and imaginary roots (1920, 148).

More generally, Leibniz's enormous capacity to create new mathematics can be attributed to his commitment to the principle of continuity, i.e., his commitment to finding similarity in difference.

> … heterogeneity for Leibniz need not threaten the rational intelligibility of things, because he can make use of the principle of continuity (a more explicitly "algebraic" version of the principle of sufficient reason) to hold heterogeneous terms together. Thus we see why a characteristic of "blind," symbolic notation should play such a central role in Leibnizian method. Descartes never seemed to realize that the algebra he helps to devise has the great virtue of allowing the hypothesis of similarity of structure while at the same time preserving the differences between the things compared. The ability of algebraic structure to allow this sameness-with-difference is its great genius and one of the reasons it has sparked so much invention. (Grosholz and Yakira 1998, 21)

This rather abstract account of the function of conceptual continuity in Leibniz's thought presents more formally the idea that Leibniz stretched concepts in order to achieve greater systematization, that is, to achieve the greatest possible unity and extension of knowledge. Moreover, the epistemic ideal set by the principle of continuity was by no means a regrettable aberration in Leibniz's philosophy. For it explains why Leibniz's infinitesimal difference approach to the calculus trumped Newton's more traditional first and last ratio approach in the development of $18^{th}$ century mathematics and physics (cf. Bos 1980, 92-3).

Infinitesimals contribute to systematization in the same manner as imaginaries and negatives, and Leibniz was quite clear about this. For instance, in comparing his geometric quadratures with his analytic quadratures, Leibniz asserted

> One part was confined to assignable quantities which had been treated not only by Cavalieri and Fermat and Fabri, but by Gregoire, Guldin and Dettonville (Pascal) as well; while the other part was based on inassignable quantities, and it promoted Geometry to a much greater extent. (GM V, 400-1; cf. Bos 1980, 62-5)

Just as imaginary roots both unified and extended the method for solving cubics[31], likewise infinitesimals unified and extended the method for quadrature so that, e.g., quadratures of general parabolas and hyperbolas[32], could be found by the same method used for quadratures of less difficult curves. Transcendental quantities, e.g., logarithms, provide another example of the systematizing power of Leibniz's calculus. Leibniz applied the term "transcendental" to logarithms because they could not be represented by algebraic equations; they transcended the resources of algebra. By supplementing ordinary algebra with the infinitesimal differences and infinite summations of his calculus, Leibniz was able to incorporate these quantities into the body of mathematics and the study of nature. For example, Leibniz observes that the equation $y = \sqrt{2x - x^2} + \int \frac{dx}{\sqrt{2x-x^2}}$

> perfectly expresses the relation between the ordinate y and the abcissa x; from it all the properties of the cycloid can be demonstrated. In this way analytic calculus is extended to those lines that have hitherto been excluded for no other cause than that they were not believed to be subject to it. (GM V, 231; translation from Mahoney 1984, 421)

Infinitesimals share, then, the virtues that led Leibniz to admire pure fictions, such as negatives, and imaginaries. We have mostly lumped these virtues under the 'meta'-virtue, systematization,[33] which covers traditional theoretical virtues such as generality, simplicity (or elegance), and the like. And like negatives and imaginaries, infinitesimals make possible an art of invention, by stretching old concepts to cover new cases. Pure fictions are prized because they contribute to the systematization of knowledge and to progress. It is our contention that even if it had not occurred to Leibniz that infinitesimals could be defended by the syncategorematic strategy, he would have had no qualms about defending them as pure fictions. Even if they cannot be eliminated in favor of ordinary quantities, they are nonetheless *fondé en réalités* because of their contribution to our knowledge of reality. We would urge, then, that *sano sensu* – in the proper sense - fictions are employed to systematize knowledge and thereby promote progress.

---

31 Leibniz essentially created determinant theory in order to solve equations of higher degree (cf. Knobloch 2006, 114-5).

32 For k constant and n, m $\in$ N, these are $y^m=kx^n$ and $ky^mx^n$, respectively.

33 Where we have used "systematization," Leibniz would employ "universality" or "universal concepts."

**§5. Formalism.** If, as we have argued, Leibniz's use of mathematical fictions is primarily for their systematic effect, rather than simply for convenience, then his fictionalism strongly anticipates the formalist approach to mathematics, championed by Hilbert, Robinson, and others. Robinson lists three essential properties of formalism.

> (1) Infinite totalities do not exist and any purported reference to them is literally meaningless; (2) this should not prevent us from developing mathematics in the classical vein, involving the free use of infinitary concepts; and (3) although an infinitary framework such as set theory, or even only Peano number theory cannot be regarded as the ultimate foundation for mathematics, it appears that we have to accept at least a rudimentary form of logic and arithmetic as common to all mathematical reasoning. (1970, 45)

These properties do not describe Leibniz's treatment of infinitesimals exactly, but this is not surprising as twentieth century mathematics employs concepts and techniques that were not exploited in the seventeenth century. But the parallels are strong nonetheless. (1) Because they involve contradictions, Leibniz rejects infinitesimals, along with negatives and imaginaries as *objects* of mathematical science. (2) Nonetheless he develops mathematics by freely employing these concepts, their use being governed by his principle of continuity.[34] (3) Common to all mathematical reasoning is a core of mathematical knowledge that concerns objects constructed in imagination from sensible content, viz., *entia rationis*. Our discussion in §§2-4 is intended as both explication and support for these three claims.

At this point it should be apparent that Leibniz's fictions play the role of ideal elements in Hilbert's sense. After arguing that there is no evidence of the actually infinite in reality, Hilbert proposes that we conceive the infinite by the "important and fruitful method of ideal elements" (1925, 185-7)

> These ideal … elements have the advantage of making the system of connection laws as simple and perspicuous as possible. … Another example of the use of ideal elements are the

[34] Boyer and others object that the principle of continuity leads to paradoxical results; e.g., according to the principle of continuity, the limit of a sequence of rationals should be rational (1959, 256). Against this we would argue that applications of the principle of continuity are constrained by the need to preserve existing knowledge. After all, the purpose of the principle is to unify knowledge rather than revise it.

> familiar complex-imaginary magnitudes of algebra which serve to simplify theorems about the existence and number of the roots of an equation. (187).

Hilbert praises ideal elements on the same grounds that Leibniz praises mathematical fictions, viz. the overall systematization which they underwrite. Furthermore, he contrasts the ideal portions of a mathematical theory with the portion obtained “through intuitive material considerations” (192). This distinction mirrors Leibniz’s distinction between fictional entities and ideal ones. It is this parallel that explains Robinson’s conflation of “ideal” and “fictional”: Infinitesimals are ideal in Hilbert’s sense but not in Leibniz’s sense; where Hilbert uses “ideal,” Leibniz would use “fictional.” Hilbert’s ideal elements and Leibniz’s fictions perform the same function, viz., the simplification of theory. This is not to suggest that Hilbert embraced infinitesimals. He was a strong proponent of Weierstrass’s reduction of propositions involving infinitely large and infinitely small quantities to statements referring only to finite quantities (183).[35] Instead of being quantities, Hilbert’s ideal elements are infinite aggregates or sets (188). This is not surprising. As mathematicians gradually abandoned the conception of mathematics as the science of quantity in the nineteenth century, sets replaced quantities.

Hilbert’s attempt to demonstrate the consistency of arithmetic is essentially an attempt to demonstrate that adjoining infinite sets to the finite or contentual truths of arithmetic leads to a reliable theory. Ishiguro’s contention that Leibniz is a proto-Hilbert arises, no doubt, from this aspect of Hilbert’s philosophy. For the syncategorematic interpretation, if it were successful, would establish that infinitesimals are reliable as long as traditional geometry is reliable. But the shortcomings of the syncategorematic interpretation, as well as the shortcomings of Hilbert’s program (i.e., Gödel’s proof that there can be no finitistic demonstration of the consistency of arithmetic) indicate that this is not the most charitable interpretation of Leibniz’s fictionalism.

---

[35] Proponents of nonstandard analysis are explicit about employing the method of ideal elements: “The study of nonstandard analysis locates itself, mathematically, in the context of the method of ideal elements. This is a time-honored and significant mathematical idea. One simplifies the theory of certain mathematical objects by assuming the existence of additional "ideal" objects as well. Examples are the embedding of algebraic integers in ideals, the construction of the complex number system, and the introduction of points at infinity in projective geometry.” (Davis 1977, 1)

It is better to conceive and evaluate ideal elements primarily in the light of contribution to mathematics, a position not so far from Costabel's suggestion that fictions are justified a posteriori. Thus we disagree with Ishiguro's 'proto-Hilbert' classification. Leibniz's best thinking about mathematical fictions is closer to Robinson's idea that infinitary concepts are literally meaningless, in the sense that they ultimately lack referents.

**Conclusion.** Infinite sets have the property that they can be shown to be equinumerous with a proper part of themselves. This contradicts the principle that the whole is greater than the part, a principle to which Leibniz was firmly committed. Infinite aggregates would, therefore, be fictions from a Leibnizian perspective. Robinson's creation of a consistent theory of infinitesimals consists, roughly, in defining the hyper-reals by set theoretic means and showing that any first order statement that holds for a real function holds for the hyper-real natural extension of that function. Thus a real function and its hyper-real extension give the same value when applied to the same real number (Keisler 1986, 28-9). To be sure, the demonstration of this result requires resources unavailable to Leibniz. But for Robinson those resources occupy the same role as infinitesimals and imaginary roots did for Leibniz. They are formal or symbolic techniques – constructed in analogy with concrete or 'ordinary' mathematics – that impart greater simplicity and systematicity to ordinary mathematics.